\documentclass[12pt,nofootinbib,aps,pra,reprint,superscriptaddress,amsmath,amsfonts,amssymb]{revtex4-1}
\usepackage[utf8x]{inputenc}
\usepackage{graphicx}
\usepackage{multirow,array}
\usepackage{amssymb}
\usepackage{float}
\usepackage{setspace}
\usepackage{amsmath}
\usepackage{textcomp}
\usepackage{url}
\usepackage[colorlinks=true, allcolors=blue]{hyperref}
\usepackage{color}
\usepackage{wasysym} 
\usepackage{MnSymbol}
\usepackage{braket}
\usepackage[all,cmtip]{xy}
\usepackage{dsfont}
\usepackage{ucs}

\definecolor{Blue}{rgb}{0,0,0.9}
\definecolor{Red}{rgb}{0.9,0,0}
\newtheorem{defini}{Definition}

\newtheorem{axioma}{Axiom}

\newcommand{\arcangleru}{%
	\mathord{<\mspace{-10.5mu}\mathrel{\rcurvearrowup}\mspace{2mu}}%
}

\newcommand{\arcangled}{%
	\mathord{<\mspace{-10.5mu}\mathrel{\lcurvearrowdown}\mspace{2mu}}%
}

\newcommand{\criss}[1]{\textcolor{orange}{#1}}
\newcommand{\bfr}[1]{\textcolor{Red}{#1}}
\newcommand{\lmg}[1]{\textcolor{magenta}{#1}}
\newcommand{\lmgg}[1]{\textcolor{teal}{#1}}
\usepackage{soul}

\begin{document}

\title{Operational approach to the topological structure of the physical space}

\author{B. F. Rizzuti}
\email{brunorizzuti@ice.ufjf.br}

\author{L. M. Gaio}

\affiliation{Depto. de F\'isica, ICE, Universidade Federal de Juiz
de Fora, MG, Brazil}

\author{C. Duarte}
\email{crsilva@chapman.edu}
\affiliation{Schmid College of Science and Technology, Chapman University, 
  One University Drive, Orange, CA 92866}

\begin{abstract}
Abstract axiomatic formulation of mathematical structures are extensively used to describe our physical world. We take here the reverse way. By making basic assumptions as starting point, we reconstruct some features of both geometry and topology in a fully operational manner. Curiously enough, primitive concepts such as points, spaces, straight lines, planes are all defined within  our formalism. Our construction breaks down with the usual literature, as our axioms have deep connection with nature. Besides that, we hope this operational approach could also be of pedagogical interest. 
\end{abstract}

\maketitle {\bf Keywords:} Operationalism, Hausdorff Spaces, Normed Spaces. Operational constructions in Mathematics. 

\section{Introduction}

Take a minute and try to come up with a definition of a ``vector". Now take a second minute and define, say, ``angle". Primitive concepts coming from geometry are widely used and taken for granted. Indeed, any educated guest has an intuition of what points, straight lines, planes and space are. Basic associated structures such as vectors and angles are on the same footing, although if you ask an undergraduate student what a vector is, the answer would likely bear the words:
``well... well... maybe something that has a direction and magnitude?''. The most experienced researcher will for sure go further and by taking a rather abstract point of view they might say either that "a vector is an element of the underlying set that forms a vector space" or that ``vectors are elements of the tangent space $T_{p}M$ on the point $p$ of a manifold $M$" or assuming a more pragmatical approach they might say something along the lines ``vectors are objects that transform according to a certain law". In either cases neither definitions are fully satisfactory, as they are completely uncorrelated with our daily-basis experience. 

Both basic and advanced literature on geometry \cite{Butkov68,Carmo2016,Spivak1999,Arfken85} does not fulfill this gap and the reader is tacitly expected to assume basic knowledge on the subject. So that, our prime question in this paper is: can we bring all ill-defined geometric primitive concepts and related ones to a conscious level? In other words, is it possible to start geometry by defining its basic ingredients? The aim of this work is to properly answer these and other related questions.

We adopt here what is known in the literature~\cite{Opera09} as operational formalism. It was extensively advocated by P. W. Bridgman \cite{Bridgman}. Roughly speaking, the idea is that no concept can be set unless one provides an experimental prescription on how to define it. Although we comprehend the tension of accepting this philosophical current, we cannot underestimate its pedagogical power.

Our program starts, then, by basic geometry, passing through the concepts of linear algebra and culminating with what we shall call ``compass-ball'' based topology. We will give not only an operational flavor to \textit{a priori} unplugged-from-the-world math structures~\cite{Leifer2016}, but also we will conclude what is the very basic topological structure our space possesses.

\begin{figure}
    \centering
    \includegraphics[scale=0.25]{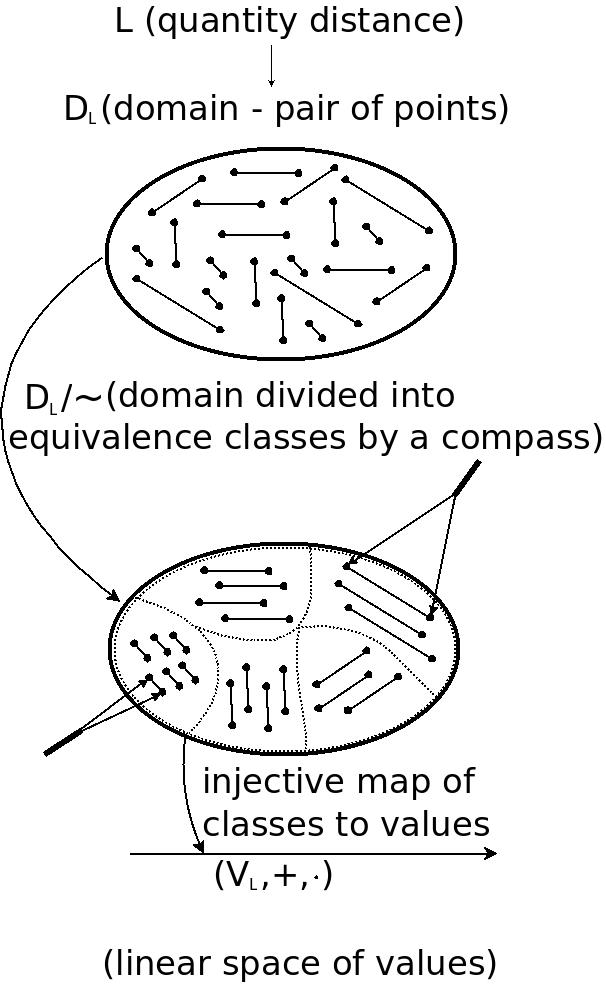}
    \caption{Schematic construction of the physical quantity distance.}
    \label{g1d}
\end{figure}

Curiously enough, the goal of our paper is intrinsically connected to the notion of physical quantities~\cite{grb,ldb2019,lesche12}. Take for instance what usually goes by the name of ``distance". In order to define such a concept, firstly we need to define the set where it makes sense, namely, the set of pair of points. We call it domain of the quantity. Secondly we do need to separate the domain by equivalent pairs of elements. For the case of distance, it is done by a compass. Finally, one associates each class to a set of values of the quantity, after all, predictions, comparisons, precise measure of sizes and so on are expected in quantitative sciences. The illustrative Figure \ref{g1d} summarizes these steps.

One should also notice that the operational approach we carry out in here is not restricted only to the realm of geometric and topological concepts, as it has a direct connection --for instance-- with what has been known as \emph{resource theories}~\cite{CFS16,TF17,DHW08}. Leveraging on a fully operational point of view, these kind of theories aim to give a precise, rigorous mathematical meaning for usually ill-defined commonly used terms. Within the scope of those theories it is possible to define what is a resource, what is a allowed transformation between resources, how to evaluate which resource is better for a given task... and so on. To the same extent resource theories have been used to clarify obscure concepts in foundations of physics~\cite{Vicente14,GA15,DA18} and in thermodynamics~\cite{TF17,BHORS13,BG15}, we expect that our formalism can also be useful to shed a new light on old and not-so-well understood math concepts.

The paper is organized as follows: Sec.~\ref{SecBasicIngredients} lays down the basic ingredients necessary for our operational approach, namely points, straight lines, and planes. In Sec.~\ref{SecDisplacementVectors} also leveraging on the operational formalism we introduce and discuss the concept of displacement vectors, and before moving on to the topology one could obtain from our framework~\ref{SecInduced}, we dedicate Sec.~\ref{SecAngles} to introduce the angle between two vectors. Finally, wrapping up our work, in Sec.~\ref{SecConclusions} we present our conclusions and discussions.

\section{Basic ingredients: points, straight lines and planes}\label{SecBasicIngredients}

Centering our attention on a radical and well-justified approach, we dedicate this section to introduce and deal with familiar, deeply rooted concepts that will be necessary in all subsequent sections. Following what part of the authors have done in Refs.~\cite{grb,ldb2019}, in here getting rid of the usual methodology and assuming an operationalist posture, we will discuss how it is possible to define objects like points, spaces, straight lines and planes, for instance, through a fully operational point of view.  

The point of view that we adopt here is that of a craftsmen or a locksmith that have in their possession a compass, a pencil and a piece of wood to work on. The piece of wood, as we will discuss below, is to be seen as a body, where by using a pencil they mark crosses on. Summing up, the idea here consists in giving operational meaning to primitive concepts like \emph{crosses, bodies...}.


\subsection{Points, rigid bodies and spaces}

We begin with the following 
\begin{defini} \label{point}
    We will call \emph{point} the center of a cross made with a pencil on a body.
\end{defini}

Leaving aside craftsmen and locksmiths, even without noticing it, we do make use of this idea. Think of hanging a painting on a wall. We have been taught that before drilling into the wall, we should first draw a little cross, determine its ``center'' and only then use the driller. Remarkably, the reader should notice that though ingenuous, the definition above brings the idea of constructing well-known mathematical objects from a operational point of view.

Now, with Def.~\ref{point} in hands, there is an idealizing assumption that will be summarized as

\begin{axioma}
    There are bodies in nature obeying the following rule: any two arbitrarily marked points on them can always be reached out by both needles of a compass keeping the same aperture constant regardless of time   
\end{axioma}
\begin{defini}
    These bodies will be called rigid bodies of reference, or simply \emph{rigid bodies}.
\end{defini}     
     
\begin{figure}
    \centering
    \includegraphics[scale=0.35]{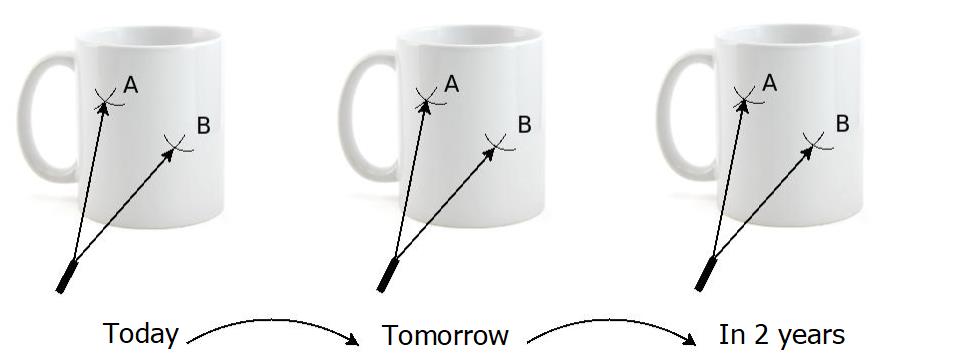}
    \caption{Distances between points in a rigid body do not change with time.}
    \label{rigbody}
\end{figure}

The sequence of images in Fig.~\ref{rigbody} illustrates this axiom. As a counter-example though, we could go with a well-known thermodynamic phenomenon, that of thermal dilation~\cite{Callen85}. Suppose we had marked points $A$ and $B$ on a golden ring during a cold winter day using a wooden compass. Over the summer, during a hot day, the compass fixed aperture would no longer fit on the marked points. Under these conditions, our golden ring would not work as a rigid body.

Imagine we are now able to unite rigid bodies with a strong enough adhesive, so that they cannot move. Take bricks for example. They are rigid bodies by excellence, and by uniting them with cement, to form a wall, the wall itself will become a new rigid body. This union shall be called a \emph{steady union}. Remarkably, notice that the steady unions can be used as a mechanism to increase the number of points our craftsmen could deal with. The more bodies we steadily put together, the more points\footnote{If your compass is such that it cannot differentiate a given marked point from another, by all means you should regard these two different points as being equal. Physically, if your measurement apparatus cannot differentiate between two objects, one should consider them as being the same.} we will end up with. 

Let us now consider the set of all rigid bodies.
\begin{axioma}
The steady union of rigid bodies is an equivalence relation.\footnote{To the reader who is not familiar with the concept of equivalence classes and relations, we recommend \cite{halmos}, which is sufficiently formal. For a more intuitive interpretation, see \cite{grb}.}
\end{axioma}

We will denominate the classes defined by the steady union of rigid bodies as a \textit{rigid frame of reference} or \textit{rigid system of reference}. We use system of reference (SR) for short. 

Going a step further, let us consider the set of points that can be marked on a frame of reference. Use your abstraction to imagine this unending quantity of points as if their underlying rigid body had been removed. With that, we gain
\begin{defini} \label{fraoref}
    The set of points in a frame of reference is called \emph{space}. We denote the space built from a certain system(frame) of reference $SR$ as $\mathcal{E}_{SR}$.
\end{defini}

Here we can already notice a rupture with the Newtonian notion of a absolute space. As a matter of fact, two distinct systems (frames) of reference define distinct spaces. There are still two subsets of points in $\mathcal{E}_{SR}$ that are of interest to us, and shall therefore be defined, namely straight lines and planes.

\subsection{Straight lines and planes}

The concept of a straight line is intimately connected to that of the physical quantity ``distance" ~\cite{Carmo2016,Spivak1999}. The details can be found in subsection 2.1 of \cite{grb} and in Ref.~\cite{ldb2019}: there the authors defined the domain $D_L$ of of the physical quantity length $L$, which are in fact pairs of points. The equivalence relation that divides $D_L$ in pairs of equivalent distances is defined by the use of a compass: a pair of points, say $(A,B)$, defines the same distance as the pair $(C,D)$ when the compass fits in both pairs, without changing its aperture. Associating the resulting equivalence classes with numerical values lying in $V_L$, they have also defined a notion of sum and multiplication by a scalar in $V_L$ (see Fig.\ref{g1d}), obtaining at the very end a structure of a vector space $(V_L,+,\cdot)$. Elements of $V_L$ are denoted by $d(A,B)$, and as we will see might be interpreted as being the \emph{distance} between $A$ and $B$. 

\begin{figure}
    \centering
    \includegraphics[scale = 0.3]{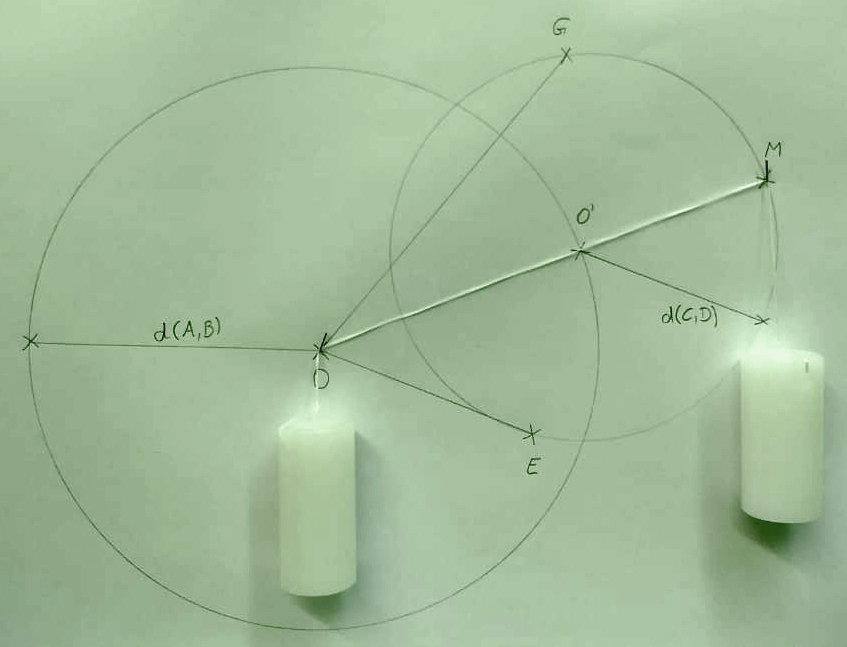}
    \caption{Illustration of the construction of (i) sum of distances, (ii) triangle inequality, (iii) straight segment and (iv) alignment of points.}
    \label{orderrel}
\end{figure}

For the sake of consistency, we shall repeat the construction of the sum operation 
\begin{equation}
+ : V_L \times V_L \rightarrow V_L.
\end{equation} 
For doing so, we start by granting $V_L$ with an order relation. Given the points $A$, $B$ and $O$ in $\mathcal{E}_{SR}$, with an open compass defined by the distance $d(A,B)$ we can form the set 
\begin{eqnarray}\label{Eq.DefBall}
    \mathcal{B}(O,d(A,B)) := \{X \in \mathcal{E}_{SR}\,|\, d(O,X) = d(A,B)\}. 
\end{eqnarray}
$\mathcal{B}(O, d(A,B))$ is called sphere of center $O$ and radius $d(A,B)$. Now, given an arbitrary point $H \notin \mathcal{B}(O, d(A,B))$, we put the compass needle point over $O$, and without taking it out from the rigid body, we draw a line connecting $O$ and $H$. In the event where this is possible without the line crossing $\mathcal{B}(O, d(A,B))$, then we say that $d(O,H) < d(A,B)$. Otherwise, $d(O,H) > d(A,B)$. Fig.~\ref{orderrel} illustrates two instances of this order relation. 
In this case, for the pair of points $E,G \notin \mathcal{B}(O,d(A,B))$ our operational construction says that $d(A,B) > d(O,E)$ and $d(A,B) < d(O,G)$. Remarkably, since the distances defined by two arbitrary pair of points may be compared by drawing the proper sphere, this order relation is a total order.

Coming back, now, to the notion of sum, to obtain $d(A,B) + d(C,D)$, we draw the set $\mathcal{B}(O', d(C,D))$, where $O' \in \mathcal{B}(O, d(A,B))$ and we search for $M \in \mathcal{B}(O', d(C,D))$ such that the distance $d(O,M)$ is maximal. Hence, we define
\begin{eqnarray}
    d(A,B) + d(C,D) := d(O,M).
\end{eqnarray}

With that in hands we naturally gain the triangle inequality, see Fig.~\ref{orderrel}
\begin{eqnarray}\label{trine}
    d(O,G) + d(G,M) \geqslant d(O,M),
\end{eqnarray}
that is valid for three arbitrary points in $\mathcal{E}_{SR}$. 

What  we are interested in though, is the equality, whose geometric interpretation lies in Def.~\ref{fraoref}. The most intriguing fact however, is that by placing an in-extensible line through $O$ and $M$ and keeping it stretched, it will go through $O'$, as shown in Fig.~\ref{orderrel}. Hence we say that $O$, $O'$ and $M$ are aligned. Moreover, we say that $O'$ is in between $O$ and $M$. This discussion motivates the following:

\begin{defini}\label{Def.LineSegment}
The \emph{line segment} between two any points $A,B \in \mathcal{E}_{SR}$ is the subset
\begin{equation}
    r(A \leftrightarrow B):=\{ X \in \mathcal{E}_{SR} | \,\,  X \,\, \mbox{is in between} \,\, A \,\, \mbox{and} \,\, B\}.
\end{equation}
In a shorthand notation, it would be useful to write down $r(AXB)$ to mean that $X$ is a point lying in between $A$ and $B$.
\end{defini}

Finally, we can define what we mean by a straight line:

\begin{defini}\label{Def.StraightLine}
Given $A \neq B \in \mathcal{E}_{SR}$, we define the straight line passing through $A$ and $B$ as the following union:
\begin{equation}\label{Eq.DefStraightLine}
    r_{AB}:= \bigcup_{X \in \mathcal{E}_{SR}} \left[ r(XAB) \cup r(AXB) \cup r(ABX) \right].
\end{equation}
\end{defini}

Essentially, the operational construction of a straight line, or even of a determined line segment, is executed by stretching an in-extensible line through two arbitrary points, which can also be seen \footnote{We used a piece of dental floss as an in-extensible thread to represent the line segment, tied to two plastic weights.} in Fig.~\ref{orderrel}.

Once a straight line is defined, before exploring what we will understand as a plane, we still need one more definition, namely, of distance between point and straight line.
\begin{defini}
    Let $P \in \mathcal{E}_{SR}$ be a point and $r \subset \mathcal{E}_{SR}$ a straight line.  We say that \begin{equation}d(P,r) := \inf \{  d(X,P)| X \in r\} \end{equation} is the distance from the point $P$ to the straight line $r$.
\end{defini}
With this definition in hands, and in complete analogy to the spheres built previously, we define one more subset of $\mathcal{E}_{SR}$.
\begin{defini}
    Let $d(A,B)$ be a distance between two points in the space and $r \subset \mathcal{E}_{SR}$ a straight line. We then define the set $\mathcal{C}(r, d(A,B)) := \{ X|d(X,r) = d(A,B)\}$ and denominate it as the \emph{cylinder} originated from $r$ and with radius $d(A,B)$.
\end{defini}
From this definition we have the following verified fact, expressed as an
\begin{axioma} \label{cylinder}
    Given an arbitrary point $P$ in $\mathcal{C}(r,d(A,B))$, there exists one unique straight line $s \subset \mathcal{C}(r,d(A,B))$ and containing $P$ such that $d(R,s) = d(A,B), \forall R \in r$.
\end{axioma}
\begin{figure}
    \centering
    \includegraphics[scale = 0.5]{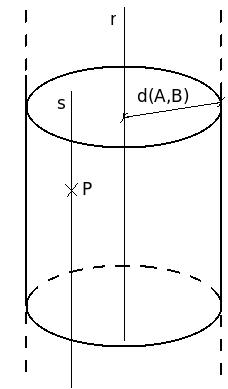}
    \caption{Representation of Axiom~\ref{cylinder}.}
    \label{cylfig}
\end{figure}
The geometric interpretation of this axiom can be seen in Fig.~\ref{cylfig}. As a consequence, we note that $d(R,s) = d(A,B), \forall R \in r$ and also $d(S,r) = d(A,B), \forall S \in s$. In this case, we say that $r$ and $s$ are \emph{parallel}. 

Finally, the last definition we need to describe our operational approach to plane geometry is that of a plane.
\begin{defini}
    Let $r, s \subset \mathcal{E}_{SR}$ be two straight lines and $P \in s$ an arbitrary point of $s$. We name the set $\Pi = r_{XP} \cup s, \forall X \in r$, where $r_{XP}$ is the straight line through $P$ and $X \in r$, the \emph{plane} defined by the straight lines $r$ and $s$.
\end{defini}
\begin{figure}
    \centering
    \includegraphics[scale = 0.35]{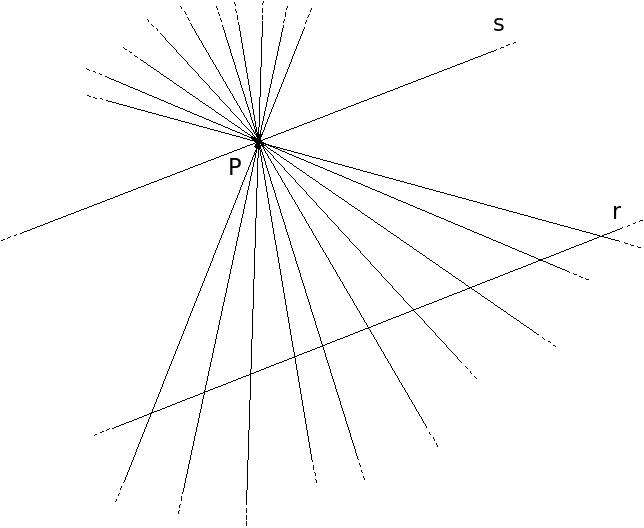}
    \caption{Operational definition of a plane.}
    \label{plane}
\end{figure}
Figure \ref{plane} illustrates this construction.

With these definitions we complete  our structuring of space, and can now proceed to the characterization of associated geometric features as well as physical quantities~\cite{grb,ldb2019}  originated from this initial step. 

\section{Displacement vectors}
\label{SecDisplacementVectors}

This Section is devoted to describe the physical quantity called displacement vectors, or, for short, vectors.  

As usual, our first step consists into characterizing its domain: it is given by ordered points of the space
\begin{eqnarray}
     \mathcal{P}_{(\cdot,\cdot)} = \mathcal{E}_{SR}\times\mathcal{E}_{SR}
\end{eqnarray}
whose elements shall be denoted $\overrightarrow{(A,B)}$.
The arrow indicates that the pair $\overrightarrow{(B,A)}$ is  different from $\overrightarrow{(A,B)}$. The name displacement vector comes from one of its possible physical interpretations. We could imagine, for instance, a particle moving from $A$ to $B$ in a straight line. In this sense, $r_{AB}$ is called the support line of $\overrightarrow{(A,B)}$. The points $A$ and $B$ will be called \emph{origin} and \emph{end} of the corresponding pair. 

\begin{figure}[h!] 
    \centering
     \includegraphics[scale=0.35]{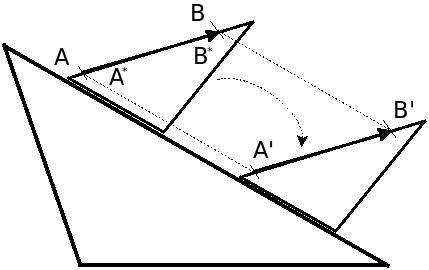}
    \caption{Parallel transport of an ordered pair of points.}
    \label{transp_p}
\end{figure}

\begin{figure}[!h]
    \centering
    \includegraphics[scale=0.4]{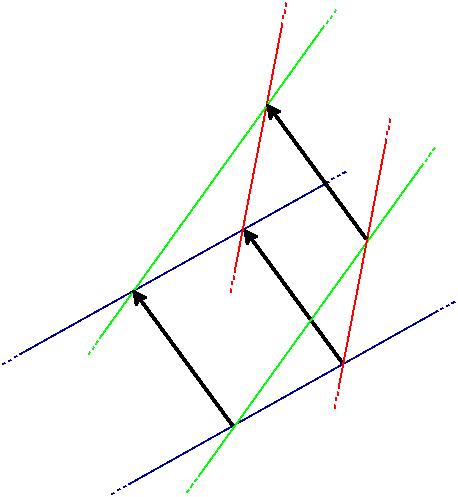}
    \caption{Parallel transport transitivity.}
    \label{transit}
\end{figure}

Moving on, our next step consists of  breaking down the domain $\mathcal{P}_{(\cdot,\cdot)}$ into classes. It is done by an experimental procedure called \emph{parallel transport}~\cite{lesche12}, denoted by
\begin{eqnarray}
     \top \subset \mathcal{P}_{(\cdot,\cdot)} \times \mathcal{P}_{(\cdot,\cdot)}.
\end{eqnarray}
Let us make a pause, though, and give  an operational prescription of how to perform it. The transport of a pair  $\overrightarrow{(A,B)}$ is done with the assistance of two set squares. We set one of them over the line segment defined by $A$ and $B$, then we mark points $A^*$  and $B^*$ on the square juxtaposed with $A$ and $B$ and slide the marked square over the other through the desired direction, see Fig.~\ref{transp_p}. 
Experimentally we may see that this in indeed an equivalence relation, for 
    
    i) $\overrightarrow{(A,B)} \top \overrightarrow{(A,B)}$, just do not move the square.
    
    ii) $\overrightarrow{(A,B)}\top \overrightarrow{(C,D)} \Rightarrow \overrightarrow{(C,D)} \top \overrightarrow{(A,B)}$, since we  could, in principle, slide the square back an forth, inverting when necessary the initial direction of the first displacement. 
    
    iii) $\overrightarrow{(A,B)} \top \overrightarrow{(C,D)}$ and     $\overrightarrow{(C,D)} \top \overrightarrow{(E,F)} \Rightarrow \overrightarrow{(A,B)} \top \overrightarrow{(E,F)}$. This transitivity may be seen in Fig. \ref{transit}.

Adopting the standard notation, we denote $\vec{a},\vec{b}, ..., \vec{c}$ each equivalence class of pairs of ordered points that may be connected by a parallel transport. We then define the set 
\begin{equation}
\mathds{D} := \mathcal{P}_{(\cdot,\cdot)}/\top,
\label{Eq.OrganizingVectors}
\end{equation}
in order to ``organize'' the set of ordered pairs of points of space. 

\begin{figure}
\begin{center}
	\includegraphics[scale=0.4]{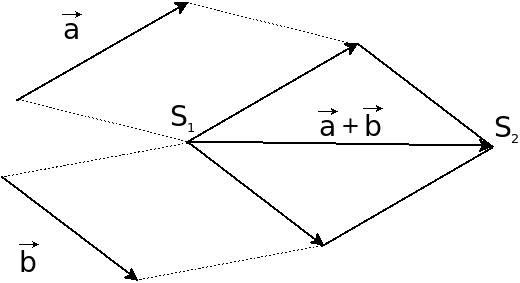}
	\caption{Constructing the vector sum.}
	\label{4}
\end{center}
\end{figure}

Following the prescription of characterizing the physical quantity we called vectors, we now define a sum of classes, 
\begin{eqnarray}
+: \mathds{D} \times \mathds{D} \rightarrow \mathds{D}.
\end{eqnarray}
Given $\vec a$, $\vec b$ and an arbitrary point $S_1$, first we connect the origins of $\vec a$ and $\vec b$ by a parallel transport to the point $S_1$. Then, we transport the pair along the support line of each one. The fact is that the end points meet in one, and only one point, say, $S_2$. The sum $\vec a + \vec b$ is denoted by $\overrightarrow{(S_1,S_2)}$. We refer to Fig. \ref{4} for the geometric picture of this construction. 
The operational construction of parallel transport together with the way we have taken $\mathds{D}$ can be used to show that the sum defined above is indeed a well-defined sum of classes. 

\begin{figure}
	\begin{center}
		\includegraphics[scale=0.4]{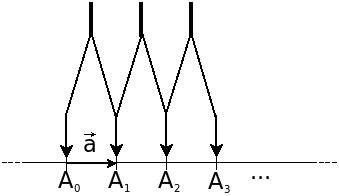}
		\caption{Constructing the multiplication of vectors by numbers.}
		\label{5}
	\end{center}
\end{figure}

One also defines the multiplication of vectors by real numbers
\begin{eqnarray}
\cdot : \mathbb{R} \times \mathds{D} \longrightarrow \mathds{D}.
\end{eqnarray}
We begin by multiplying a vector $\vec{a}$ by a natural number. Along the support line of $\vec{a}$, we mark with a compass consecutive points, keeping the compass aperture fixed, defined by the ends of $\vec a = \overrightarrow{(A_0,A_1)}$, according to Fig. \ref{5}. Each point $A_ {n + 1}$ is obtained from $A_n$. Thus, we define
\begin{eqnarray}
n \vec a := \overrightarrow{(A_0, A_n)},\, \forall n \in \mathbb{N}.
\end{eqnarray}

The multiplication by integers is totally analogous. We simply reverse the direction of $\vec a$ when the integer is negative, that is, given a negative $z$ in $\mathbb{Z}$, we have $z=-n$, for some $n \in \mathbb{N}$. Then
\begin{eqnarray}
z \vec a = (-1)n \vec a = (-1) n \overrightarrow{(A_0, A_1)} := \overrightarrow{(A_n, A_0)}.
\end{eqnarray}

\begin{figure}
	\begin{center}
		\includegraphics[scale=0.45]{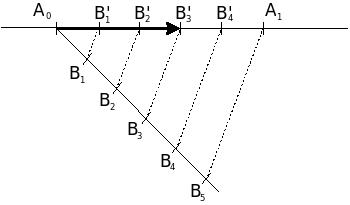}
		\caption{Geometric representation of the multiplication of a vector by a rational number.}
		\label{6}
	\end{center}
\end{figure}

We now turn to the multiplication by rationals $\frac{m}{n} \in \mathbb{Q}$, with $m\leq n$. Starting from $A_0$, we draw a line segment $r$ in any direction other than the line that supports $\vec a = \overrightarrow{(A_0, A_1)}$. With a compass opened arbitrarily, we mark the points $B_1$, $B_2$,..., $B_n$ over $r$, so that $d(A_0, B_1)= d (B_1, B_2) = ... = d(B_{n-1}, B_n)$. We draw the segment connecting $B_n$ to $ A_1 $ and trace segments parallel to the segment defined by $B_n$ and $A_1$, ranging from $B_i $, $i = 1, 2 $, ..., $n -1$ until you touch $\vec a$, where we set the points $B'_i$. Thus,
\begin{eqnarray}
\frac{m}{n}\vec a :=\overrightarrow{(A_0, B'_m)}. 
\end{eqnarray}  
The geometric representation for such a definition (multiplication of $\vec a$ for $3/5$) can be seen in Figure \ref{6}. 
The generalization for $m>n $ is straightforward. It suffices to write $m=n+m'$, so that $m/n=1+m'/n$, with $m'<n$. So, now we just apply the reasoning discussed above.

The multiplication by reals in $\mathds{D}$ may be generalized invoking the density of $\mathbb{Q}$ in $\mathbb{R}$~\cite{Rudin87}. However, we observe an experimental limitation for constructing a distance given by an irrational value. In effect, one measures the diagonal of a square formed by sides of $1m$ long. What is found? $\sqrt{2}m$? Or $ 1.41m$? A more detailed discussion of uncertainties can be seen in \cite{grb}. 

Finally, we conclude this section with a remarkably operational fact: having called the ordered pairs of vectors previously was not a mere coincidence,  the physical space $\mathcal{E}_{SR}$ associated with a given rigid reference body does have the properties of a structure well known to mathematicians -- namely, an \emph{affine space}.\footnote{We find the same structure in relativity theory. Minkowski's space-time $\mathbb{M}$ has the structure of a four-dimensional affine space \cite{oneil}.} It is defined as follows. Let $E$ be a set whose elements are called points and $V$ a vector space. Also consider an application $\vec{\cdot}: E \times E \rightarrow V$, $(p,q) \mapsto \vec{pq}$. We say that $E$ is an $n$-dimensional affine space associated with the vector space $V$ when it bears the following properties

\begin{itemize}
    \item [(i)] Given a point $p \in E$ and a vector $\vec v \in V$, there exists one, and only one point $q \in E$ such that $\vec{pq}= \vec v$.
    
    \item [(ii)] For any $p,q,r \in E$ one has $\vec{pq}+\vec{qr}=\vec{pr}$.
\end{itemize}

This reasoning makes clear the existence of a bijection between $\mathcal{E}_{SR}$ and $\mathds{D}$. In fact, given two points $A$ and $B$ in $\mathcal{E}_{SR}$, they define the corresponding  class of displacement vectors. Conversely, given a point $A \in \mathcal{E}_{SR}$ and a vector $\vec a \in \mathds{D}$, we parallel transport the origin of $\vec a$ to the point $A$, and we look at the unique end of the element in the class $\vec a$, namely, $\overrightarrow{AB}$. It is customary to write such a result as
\begin{eqnarray}
B:= A + \vec{a}.
\end{eqnarray}
Although $\mathcal{E}_{SR}$ and $\mathds{D}$ are in an one-to-one co, they are not the same. We hope the discussion above makes it more explicit.

\section{Angles}\label{SecAngles}

Before continuing to the topology of space we must give a step back though and describe a one-dimensional physical quantity that will be needed ahead: angles. For that we shall give heuristic arguments, followed by operational procedures, which are mathematically formal. 

\subsection{Constructing the physical quantity angle}

Consider the usual act of opening the fridge door to take out a bottle of water. In order to spare your own energy, you open the door just enough to take out the bottle. If you wanted to take a pan however, you would have to open the door a little wider than before. Although this example seems to be quite natural, hold on a second and think about it, what  does one mean by wider openings? Moreover, even without having defined opening, we created a order relation by saying ``wider than before". Our objective will be therefore, to formalize this notion of opening, constructing the physical quantity we will call \emph{angle}.

With the initial example, it becomes clear that the angle's domain are pairs of vectors. By restricting the initial definition to a two dimensional vector space $\mathds{D}|_{\Pi}$ that has $\Pi$ as affine space, that means, for all 
$\vec a \in \mathds{D}|_{\Pi}$ and $P \in \Pi$,  we have $P+\vec{a} \in \Pi$.  By the operational definition of plane, it becomes clear that two vectors define a plane. Therefore, we write the domain as
\begin{eqnarray}
D_A=\{ (\vec a,\vec b); \vec a, \, \vec b \in \mathds{D}|_\Pi \}. 
\end{eqnarray}

\begin{figure} [!h]
	\centering
	\includegraphics[scale=0.45]{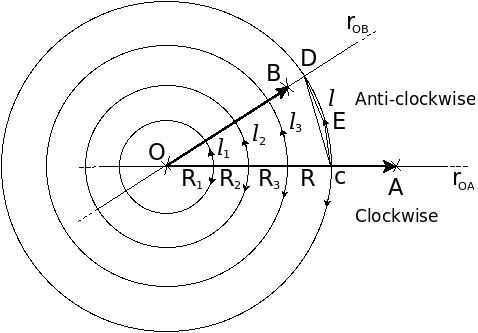}
	\caption{Operational construction of an angle.}
		\label{ang_ast}
\end{figure}

Now, we will follow a series of operational steps to characterize this new quantity. Given $\vec a = \overrightarrow{(O,A)}$ and $\vec b = \overrightarrow{(O,B)}$,

i) We trace the lines $r_{OA}$ e $r_{OB}$;

ii) We put the needle of a compass on $O$ and we draw the arc that connects the line  $r_{OA}$ to the line $r_{OB}$. This can be done in only two ways, as indicated by Fig.\ref{ang_ast}, namely,  clockwise or anti-clockwise; 

iii) Keeping the needle point on $O$, we then move the pencil lead away drawing more arcs as indicated above in step (ii), choosing one of the directions.

Surprisingly, the radiuses $R_1$, $R_2$,..., $R$ of the different arcs in Fig.~\ref{ang_ast} are different, as are their lengths as well. However, the opening defined by the vector pair remains the same. Moreover, it is an experimental fact that 
\begin{eqnarray}\label{an00}
\frac{l_1}{R_1} = \frac{l_2}{R_2} = \frac{l_3}{R_3} = \frac{l}{R}.
\end{eqnarray}
As the proportion between the arc length and the radius is independent of how much we open the compass, we can define
\begin{eqnarray}\label{an01}
\alpha := \arcangleru (\vec a, \vec b) = \frac{l}{R}.
\end{eqnarray}
In Eq.~\eqref{an01}, $l$ and $R$ correspond to the arc length and the radius as defined by the experimental procedures (i)-(iii) above. The symbol $\arcangleru$ indicates that the compass starts in line $r_{OA}$ and goes all the way to the line $r_{OB}$ anti-clockwise. It is a consequence of Eq. (\ref{an01}) that
\begin{eqnarray}
\arcangleru (\vec a, \vec b) = \arcangled (\vec b, \vec a). 
\end{eqnarray}
Naturally Eq. (\ref{an01}) defines the quantity's numerical value and, in addition, $\alpha \in \mathbb{R}$. Furthermore, two pairs $(\vec a,\vec b)$ and $(\vec c,\vec d)$ are \emph{angle-equivalent} when
\begin{eqnarray}\label{an02}
(\vec a,\vec b)\sim_A (\vec c,\vec d) \Leftrightarrow \arcangleru (\vec a,\vec b) = \arcangleru (\vec c,\vec d).
\end{eqnarray}
Before we go on to define the sum and multiplication by scalars operations, a brief commentary should be made. The physical quantity distance has pairs of points as its domain. Therefore, it is necessary to clarify the meaning of ``arc length". As initially we could only measure the distance between two points, we may try to define the arc length described by a compass coming out of $C$ and ending in $D$ as
\begin{eqnarray}\label{an03}
l:= d(C,D),
\end{eqnarray}
where the points $C$ and $D$ are indicated in Fig.~\ref{ang_ast}. But the arc does not coincide with the line segment defined by the pair $C$, $D$. Therefore, the prescription that led us to Eq. (\ref{an03}) is not ideal. For fixing this, we might put one more point over the arc, say $E$, and we write
\begin{eqnarray}
l:= d(C,E) + d(E, D).
\end{eqnarray}
By the triangular inequality, $d(C,D)<d(C,E)+d(E,D)$ but we still do not reach the arc with two segments defined by $C$, $E$ e $E$, $D$. Proceeding in this way, we continue marking down more and more points $C_i$, $i=1,...,N$, for any $N \in \mathbb{N}$, on the arc, as close as the compass' aperture allows us to measure the distance $d(C_i,C_{i+1})$. Therefore, by doing so, after many iterations, we write this sum of distances as
\begin{eqnarray}
S_P(N) = \sum_{i=1}^{N-1} d(C_i, C_{i+1}),
\end{eqnarray}
with the first point being $C$ and the last being $D$. 

This process of marking points on the arc is called partition, justifying the index $S_P$. As a result of the triangular inequality the more points are selected, the closer we can get from the actual arc length:
\begin{eqnarray}
S_P(N) \leq S_P(N+1).
\label{Eq.GettingCloserArcLength}
\end{eqnarray}
This inequality between sums indicates that we should take some sort of optimization to define the length we are looking for:
\begin{eqnarray}
l:= \sup_{\mbox{$P$ is partition}} S_P.
\end{eqnarray}
In this manner, we are able to reach as close as possible to what we can call an \emph{arc length}. 

\begin{figure}
    \centering
    \includegraphics[scale=0.075]{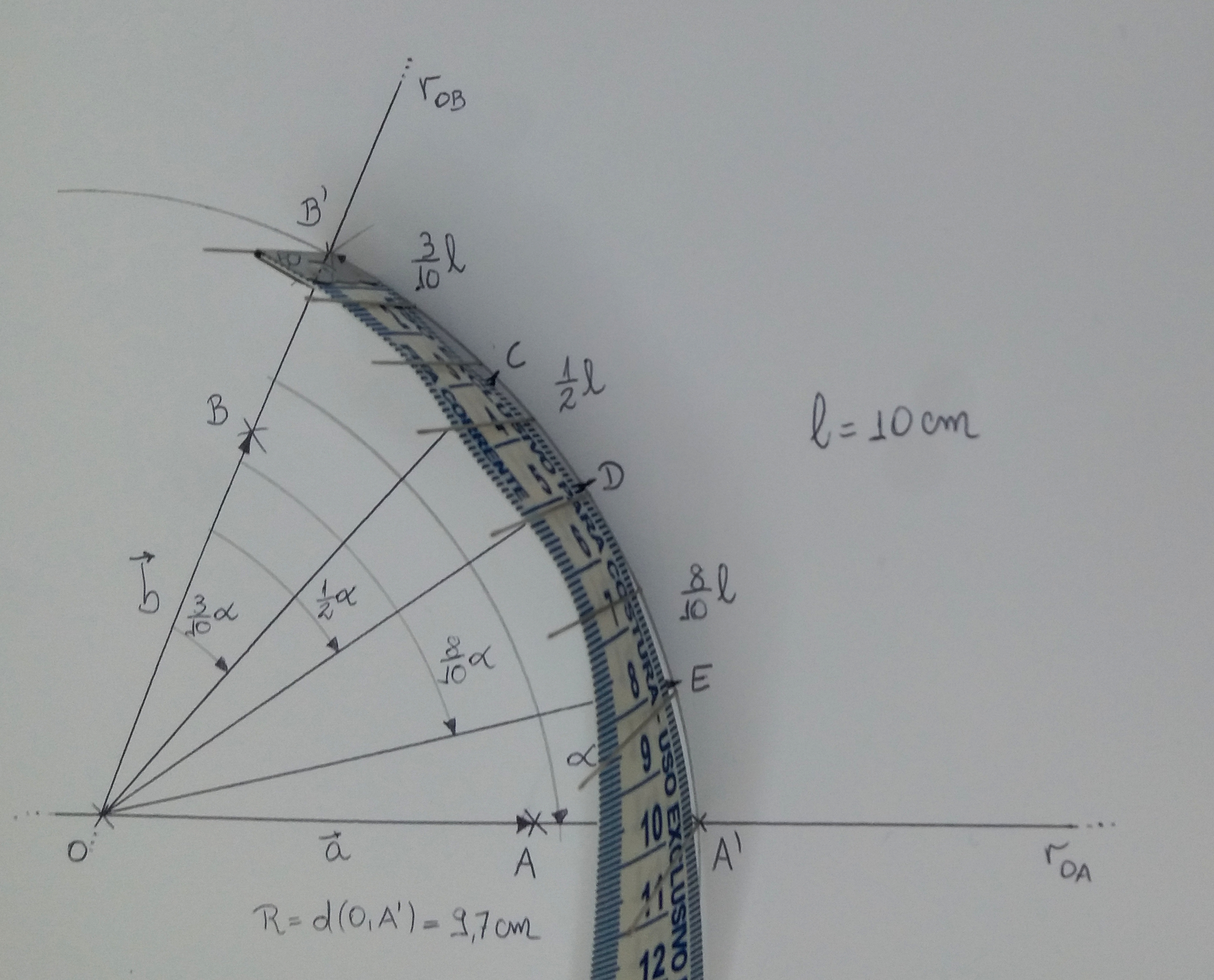}
    \caption{Constructing the multiplication by a number.}
    \label{Fig.ConstructingMultiplicationAngle}
\end{figure}

The curious practical fact is that the entire process described above results, considering the uncertainty, in the same measure we would obtain if we had used a flexible tape measure, as that of a tailor. For that, we just need to adjust the tape over the arc which we wish to measure, as framed in Fig.~\ref{Fig.ConstructingMultiplicationAngle}: the needles represent points of a possible arc partition.

Moving on, to construct the sum of angles, we still need the concept of rotation. The situation here is analogue to what has been done for the sum of distances. Writing $2cm + 3cm = 5cm$ is meaningless in the geometric sense, as given the points $A$, $B$, $C$, $D$ such that $d(A,B) = 2cm$ e $d(C,D)=3cm$, we do need a procedure to determine where the points $O$ and $M$ are, so that
\begin{eqnarray}
d(A,B)+d(C,D) = d(O,M). 
\end{eqnarray}
With the angles $\alpha = \arcangleru (\vec a, \vec b)$ and $\beta =\arcangleru (\vec c, \vec d)$ where $\vec a$, $\vec b$, $\vec c$ and $\vec d$ are in $\mathds{D}|_\Pi$, we could ask the same question: which is the pair  $(\vec e, \vec f)$ such that
\begin{eqnarray}
\arcangleru (\vec e, \vec f) = \alpha + \beta? 
\end{eqnarray}
To answer this question, we begin with a
\begin{defini}
Given $\Pi \subset \mathcal{E}_{SR}$, we call a \emph{rotation} around a point $P \in \Pi$ the mapping 
\begin{eqnarray}
R_P: \Pi &\longrightarrow& \Pi \cr
A &\mapsto & R_P(A) 
\end{eqnarray}
so that

i) $d(P,A) = d(P, R_P(A))$;

ii) $R_P(P) = P$.
\end{defini}

Mirror reflections through $P$ preserve distances, but we will exclude them from the above definition. The reason for this is that our physical space was constructed from rigid bodies of reference. It is possible to rotate a rigid body, preserving distances. But we cannot rotate our right hand to get our left hand without obtaining a new rigid body; only a reflection can do so.\footnote{Technically, the exigence of excluding reflections means asking that the mapping $R_P$ be continuous and that it can be deformed continually to the identity operator. We will discuss the operational meaning of continuity in a subsequent work.}

By identifying $\mathcal{E}_{SR}$ with $\mathds{D}$, the mapping $R_P$ induces a natural rotation of vectors.

\begin{defini}
We call a \emph{rotation} of a vector $\vec a \in \mathds{D}\vert_{\Pi}$ the following mapping, 
\begin{eqnarray}
\mathcal{R}: \mathds{D}\vert_{\Pi} &\longrightarrow& \mathds{D}\vert_{\Pi} \cr
\vec a &\mapsto & \mathcal{R}(\vec a):= \overrightarrow{(P, R_P(A))}.
\end{eqnarray}
\end{defini}

It is important to note that it is not necessary to keep one single point fixed, as we can always transport one element from the class $\vec a$ to the point $P$. From this point on, when we deal with a rotation, we will be referencing the vector rotation.\footnote{Rotations are described by theory of representations of the group $SO(3)$. The details can be seen in \cite{adb}.} 

One fact associated with rotations can be expressed by

\begin{axioma}\label{ax.rot}
Given $\vec a$ and $\vec b$ distinct classes in $\mathds{D}\vert_{\Pi}$, there exists one rotation such that $\mathcal{R}(\vec a)$ and $\vec b$ are supported by parallel lines.
\end{axioma}

Axiom \ref{ax.rot} can be visualized in an analog clock, when the seconds pointer reaches, for example, the minutes pointer.

Still using vector rotations, we can induce the rotation of a pair $(\vec a, \vec b) \in \mathds{D}\vert_{\Pi} \times \mathds{D}\vert_{\Pi}$ defining it through

\begin{eqnarray}
(\vec a, \vec b) \mapsto \mathcal{R}(\vec a, \vec b):= (\mathcal{R} \vec a, \mathcal{R} \vec b).
\end{eqnarray}
By construction, the notion of rotation we have just defined preserves angle values.
\begin{eqnarray}
\alpha= \arcangleru (\vec a, \vec b) = \arcangleru (\mathcal{R} \vec a, \mathcal{R} \vec b).
\end{eqnarray}

\begin{figure} [!h]
	\centering
	\includegraphics[scale=0.35]{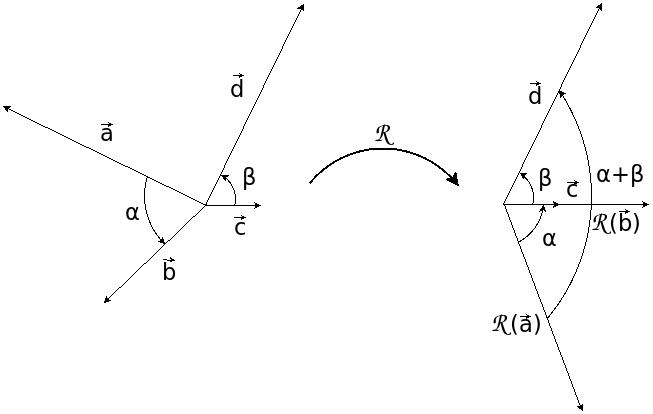}
	\caption{Geometric representation of angle sum.}
		\label{Fig_soma_ang}
\end{figure}

With this arsenal in hands, we can already characterize geometrically sums of angles. Let us take $\alpha = \arcangleru (\vec a, \vec b)$ and $\beta = \arcangleru (\vec c, \vec d)$. We search the rotation $\mathcal{R}$ so that $\mathcal{R}(\vec b)$ and $\vec c$ are parallel and with the same origin. Applying the rotation induced by $\mathcal{R}$ on the pair $(\vec a, \vec b)$, and maintaining $(\vec c, \vec d)$ fixed. We define the \emph{sum} of $\alpha$ and $\beta$ as being:
\begin{eqnarray}
\alpha+ \beta := \arcangleru (\mathcal{R} \vec a, \vec d).
\end{eqnarray}
The figure \ref{Fig_soma_ang} shows this construction. 

The \emph{multiplication} of angle values by scalars is imported from the multiplication of distances by numbers. We use once again our tailor's tape measure, with the following procedures: given an angle value $\alpha = \arcangleru (\vec{a}, \vec{b})$, we draw an arbitrary arc that connects the lines supporting the vectors $\vec{a}$ e $\vec{b}$, following the orientation defined in the construction of the angle. Measuring the arc $l$ and the radius $R$, in such a way that $\alpha = l/R$. Given $\lambda \in \mathbb{R}$, with $0<\lambda<1$, we mark on the flexible tape measure the value $\lambda\cdot l$ and return the tape over the arc, such that
\begin{eqnarray}\label{ang_n}
\lambda \cdot \alpha := \frac{\lambda \cdot l}{R}.
\end{eqnarray}
Image depicted in Fig.~\ref{Fig.ConstructingMultiplicationAngle} summarizes the procedure that results in definition (\ref{ang_n}). In this case, 
\begin{eqnarray}
\alpha = \arcangled (\vec b, \vec a) = 10cm/9,7cm = 1,03. 
\end{eqnarray}
The points $C$, $D$ e $E$ have been selected accordingly to the tape measure's scale, so that
\begin{align}
\arcangled (\vec{b}, \overrightarrow{OC})&= \frac{3}{10} \alpha \\ 
\arcangled (\vec{b}, \overrightarrow{OD}) &= \frac{1}{2} \alpha \\
\arcangled (\vec{b}, \overrightarrow{OE}) &= \frac{8}{10} \alpha.
\end{align}

We must be careful for the cases where $\lambda<0$. As we want to measure ``openings", the values of negative angles correspond to moving the compass in another direction. The analogy becomes complete when we define the space of distance values. There, we could multiply distances by real negative numbers and interpret a negative distance as, for example, a position coordinate. Other than that, thanks to the fact that
\begin{eqnarray}
\arcangleru (\vec{a}, \vec{b}) + \arcangleru (\vec{b},\vec{a}) = 2\pi,
\end{eqnarray}
we must be cautious as well when $\lambda>2\pi$, indicating that we have already turned around. This is necessary because the same pair of vectors can represent more than one angle value. To go around this lack of injectiviness between equivalence classes in $D_A/\sim_A$ defined by the relation (\ref{an02}) and values, we can add an index to the symbol $\arcangleru$
\begin{eqnarray}
\arcangleru_k (\vec{a}, \vec{b}) = \alpha + 2k\pi; \, k\in \mathbb{Z}.
\end{eqnarray}

As angle values are real numbers, we can import their usual properties, making sure that the set of values is a vector space. It is important to notice that, to distinguish a number from another which represents an opening, we use, conventionally, the word \textit{radian}. This does not mean that ``rad'' represents a physical unity. The basis of the angle value space $V_A$ is the unitary set $\{1\}$; $1\in \mathbb{R}$.

With that, we complete the description of the angle quantity. The next subsections are dedicated to some usual objects, that are extremely connected with this quantity.

\subsection{Trigonometric functions}

\begin{figure} [!h]
	\centering
	\includegraphics[scale=0.5]{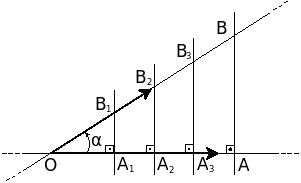}
	\caption{Operational construction of trigonometric functions.}
		\label{sencos}
\end{figure}

The idea of an angle comes from the constancy observed when we divide the arc length by the radius, expressed in (\ref{an01}). In a similar manner, we could draw line segments instead of arcs, creating, from a determined pair of vectors, triangles with angles of $\pi/2$, as indicated by Figure \ref{sencos}. Said angles are represented by a square with a dot inside, next to the points $A_1$, $A_2$, etc.
The experimental fact is that, having a fixed angle, the following operations are constant:
\begin{eqnarray}\label{seno}
\frac{d(A_1,B_1)}{d(O,B_1)}=\frac{d(A_2,B_2)}{d(O,B_2)}=\frac{d(A_3,B_3)}{d(O,B_3)} = \frac{d(A,B)}{d(O,B)}. 
\end{eqnarray}
as well as
\begin{eqnarray}\label{cosseno}
\frac{d(O,A_1)}{d(O,B_1)}= \frac{d(O,A_2)}{d(O,B_2)} =\frac{d(O,A_3)}{d(O,B_3)} = \frac{d(O,A)}{d(O,B)}.
\end{eqnarray}
With that, we defined the functions \emph{sine} and \emph{cosine}, given by the relation between the sides of straight triangles, as exposed in Eqs. (\ref{seno}) and (\ref{cosseno}) respectively. From this, the theory of trigonometric functions with an operational proposition is born.

\section{Induced mathematical structures based on operational formalism}\label{SecInduced}

Wrapping up our work we dedicate this section to the description of some important structures arising both in $\mathcal{E}_{SR}$ and in $\mathds{D}$, namely that of topology and normed/metric spaces, much in the spirit of what has been done in Refs.~\cite{grb,ldb2019}. As before, we will maintain our operational approach, giving concrete examples to illustrate the reason behind our definitions, giving therefore some physical meaning to the objects presented here.

\subsection{Compass-Based Topology}\label{SubSecCompassBasedTopology}

The introduction of a compass in the physical space is a clear operational way of defining a notion of a ball on it (see Eq.\eqref{Eq.DefBall}), and remarkably this feature already allows us to define an important set of subsets of $\mathcal{E}_{SR}$:
\begin{align}
    \tau := \{\mathcal{A} \subset \mathcal{E}_{SR} |\,\forall A \in \mathcal{A},\,\exists \, \epsilon \, \mbox{in} \, V_L, \, \mbox{with} \, \epsilon > 0; \,B(A,\epsilon) \subset \mathcal{A} \}.
\end{align}
In the above definition, by $B(A,\epsilon)$ we mean
all the points $X$ of $\mathcal{E}_{SR}$ such that $d(X,A)<\epsilon$. The pair $(\mathcal{E}_{SR},\tau)$ is called a \emph{topological space}.

\begin{figure}
    \centering
    \includegraphics[scale=0.3]{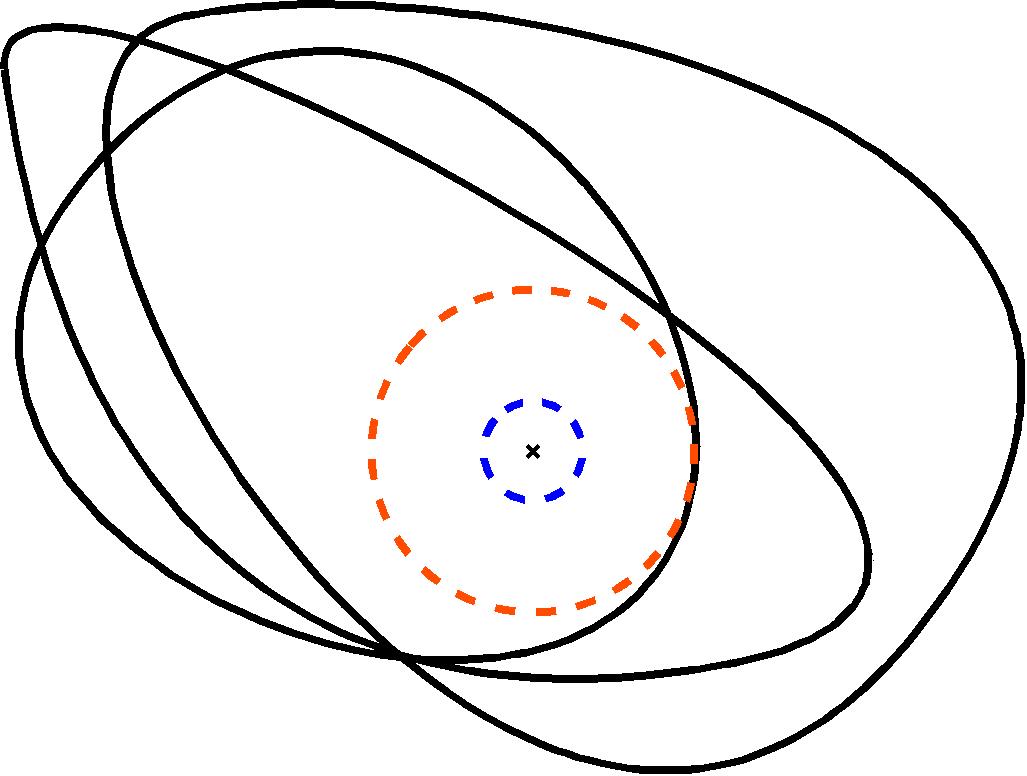}
    \caption{A point in an intersection of three open sets. The orange one would be used mainly by mathematicians, while the blue one would be used by extra careful physicists.}
    \label{open_uncert}
\end{figure}

The fact that $\tau$ really is a well-defined topology over $\mathcal{E}_{SR}$ follows the same geometric idea used in mathematical analysis to define open sets in $\mathbb{R}^n$, and as such, the proof is effectively the same\footnote{As a matter of fact we should have considered not only $\tau$ but an extension of it containing both the empty set and $\mathcal{E}_{SR}$ itself.}. Therefore, avoiding deviation from our main topic, we will not show it here. The reader may take a look at Refs.~\cite{Munkres2000} if they are interested in a proof of this fact. Nonetheless, we must make a brief commentary about uncertainties, as it is relevant during this demonstration, namely where we need to show that the topology is closed for finite intersections. Usually we take the smallest radius of all the balls centered around any point of the intersection and say that this ball, built with this radius is contained within this intersection. Mathematicians would usually be satisfied with this, as this ball would be as close to the border as we want without touching it, but a more experimental oriented physicist might be a little more cautious and, for example, divide the radius of the sphere by two, just to be completely sure that the sphere will be entirely contained within the intersection. Fig \ref{open_uncert} illustrates this difference.

Using the operationalist philosophy, we are forced to be as cautious as our imaginary physicist, when we demonstrate this fact, as there are always a certain level of uncertainty on where the border actually is. Apart from this fact, the rest of the proof follows similarly, as it is already quite operationalist in its essence. 

The topology we defined above can be used to show another interesting property of $\mathcal{E}_{SR}$, which states that given any two points in the physical space, we can find neighbourhoods for each of these points which are also disjoint from each other. When a topological space has this property it is called a \emph{Hausdorff space}~\cite{Munkres2000}.

To show that $\mathcal{E}_{SR}$ has this property, given any two distinct points $A$ and $B$ in it, we simply take two balls with  radiuses $\frac{d(A,B)}{2}$ each, and those sets are in fact disjoint (as these balls are open sets, there is no problem with the point located in the center of the line connecting $A$ and $B$, but if you want to be cautious as our imaginary physicist from before, you may always divide $d(A,B)$ by an $n \in \mathbb{N};\,n>2$).

\subsection{Metric Space}\label{SubSecMetricSpace}

Once again we are guided by the compass to insert one more specific map, which will provide $\mathcal{E}_{SR}$ with the very structure of a metric space. Let us consider the following map
\begin{eqnarray}
    d: \mathcal{E}_{SR} \times \mathcal{E}_{SR} &\rightarrow& V_L \cr 
    (A,B) &\mapsto& d(A,B)
\end{eqnarray}
This quantity clearly has the properties: 
\begin{enumerate}
    \item $d(A,A)=0$ for every point $A$ of space as the compass has no aperture.
    
    \item Any aperture of a compass defined by a pair of distinct points always defines a positive distance.
    
    \item One may use the compass to find $d(A,B)$ and changing the needles from $A$ to $B$ and vice-versa does not change the aperture, that is, $d(A,B) = d(B,A)$.
    
    \item Finally, triangle inequality also follows, as discussed before, see Eq. (\ref{trine}).
\end{enumerate}
So, we can conclude that the pair $(\mathcal{E}_{SR}, d)$ is a metric space, as claimed before. 

\subsection{Normed Space}\label{SubSecNormedSpace}

From this point on, we will be dealing with properties related to vector spaces. We will now be working with $\mathds{D}$ instead of $\mathcal{E}_{SR}$. As we stated before this is not a problem, as they are both representations of the physical space, differing only by the elements used to describe it. With this in mind we can proceed to the next step.

The essential aim of this Subsection is that of showing that $\mathds{D}$ is a normed space. As vectors are defined as classes of ordered points, say $A$ and $B$ represents $\vec{a}$, we simply define the norm of it as the distance between these two points: 
\begin{align}\label{Def.EqDefiningNorm}
    \Vert \cdot \Vert : \mathds{D} &\rightarrow V_L \nonumber \\
    \vec{a} &\mapsto \Vert\vec{a}\Vert = d(A,B).
\end{align}
We can show that this definition obeys each of the requirements to be called a norm, as given the vectors $\overrightarrow{(A,B)}$ and $\overrightarrow{(C,D)}$ in $\mathds{D}$ we can show that:

\begin{enumerate}
    \item $\Vert\overrightarrow{(A,B)}\Vert = d(A,B) \geq 0$ and $\Vert\vec a\Vert = 0 \iff \vec a = \vec 0$. This simply comes from the definition of distances  as performed in Ref.~\cite{grb}.
    
    \item By importing the triangle inequality, see Eq.~\eqref{trine}, one writes $\Vert\vec{a}+ \vec{b}\Vert \leq \Vert\vec{a}\Vert+\Vert\vec{b}\Vert$.
    
    \item $\Vert\lambda \vec{a}\Vert = \vert \lambda \vert \Vert\vec{a}\Vert$, for all $\lambda \in \mathbb{R}$. This property comes the very construction of the notion of distance, as thoroughly discussed in Ref.~\cite{grb}.  
\end{enumerate}    

The definition of norm, as we provide in Eq.~\eqref{Def.EqDefiningNorm} above, with the conditions 1-3 guarantee that the pair $(\mathds{D}, \Vert \cdot \Vert)$ is a well-formed normed space.

\section{Conclusions}\label{SecConclusions}

Fulfilling the gap normally affecting our students, in this paper we took the operationalist point of view to give a more practical, though rigorous, meaning to geometric constructs. Namely, through well-justified operational steps we came up with a robust definition of rigid bodies, spaces, straight lines, planes and finally (displacement) vectors and angles. Other than that, with the assistance of our idealized compass, we also discussed how to define a ball and therefore we ended up with a definition of a distance. Both the former and latter provided the necessary toolkit to give a step further in our formalism and see the physical space $\mathcal{E}_{SR}$ and the space of displacement vectors $\mathds{D}$ as being a Hausdorff topological space and a normed space (respectively).  

Connected with the topological structures we attached to our spaces, although we did not have approached this topic in this present work, we truly believe that we could have given one step further. Defining the inner product via the usual expression
\begin{equation}\label{Eq.DefInnerProduct}
    \vec{a} \cdot \vec{b} := \frac{\Vert \vec{a} + \vec{b}  \Vert^{2} - \Vert \vec{a} - \vec{b}  \Vert^{2}}{4}
\end{equation}
we would end up getting on the r.h.s a bilinear, symmetric, positive definite form. However, this might only be obtained if one used the ``parallelogram law". As we could not find a well-justified operational justification for the latter, we left this point to be explored elsewhere in future works.

Another fact that will be postponed for a further discussion is the connection of the inner product with the geometric structure $\mathds{D}$ possesses, that is, 
\begin{equation}
 \vec{a}\cdot \vec{b} = \Vert \vec{a} \Vert \Vert \vec{b} \Vert \cos \left[\arcangleru (\vec{a}, \vec{b}) \right].
\end{equation}
We suspect, in an operational way, that $\mathds{D}$ may be called a Hilbert space. Issues concerning the space completeness will also be addressed elsewhere. 

Also, as a matter of fact, one very natural question might have arisen right after the introduction of the space $\mathds{D}$ of displacement vectors. One could have asked about the dimension $d$ of such space. Although it is in a bijection of $\mathcal{E}_{SR}$, an affine space we usually regard as having a three-dimensional character, we do not have a better argument for this question. We do not see that as a drawback of our approach, that simply shows the limitations the operationalist model we took to describe the mathematics of the world we live in.

\section*{Acknowledgments}
This work is supported by Programa Institucional de Bolsas de Iniciação Científica - XXXI BIC/UFJF-2018/2019, project number ID45249. CD has been supported by a fellowship from the Grand Challenges Initiative at Chapman University. BFR would like to express his gratitude for the warm hospitality of the Institute for Quantum Studies and the Grand Challenges Initiative at Chapman University, where this work was concluded.


\bibliography{mybiblio}{}

\bibliographystyle{plain}

\end{document}